\documentclass[fleqn]{article}\newcommand{\ispreprint}{1}
\usepackage{amsmath,amssymb,amsfonts,latexsym,amsthm,bm,ifthen}
\usepackage{natbib}
\newcommand{\dbt}[1]{\ifcase#1{}\or.5\or1\or1.5\or2\or2.5\or3\or3.5\or4\or4.5\fi}

\renewcommand{\-}[1]{\mskip-\dbt#1mu}
\renewcommand\![1]{\bm #1}

\def\ClAuDiAeatblanks{\global\@ignoretrue}
\def\ClAuDiAstop{\relax}
\def\ClAuDiAsplit #1#2\stop{\string #1}
\def\[#1\]{\protect{\if\ClAuDiAsplit #1\stop[%
{\ClAuDiAwithlabel #1\ClAuDiAstop}\else%
{\ClAuDiAwithoutlabel #1\ClAuDiAstop}\fi}\ignorespaces}%
\def\ClAuDiAwithlabel[#1]#2\ClAuDiAstop{\protect{\begin{equation}\label{#1}\begin{split}#2\end{split}\end{equation}}}
\def\ClAuDiAwithoutlabel#1\ClAuDiAstop{\protect{\begin{equation*}\begin{split}#1\end{split}\end{equation*}}}

\newcommand{\set}[2]{\bigl\{\mskip1mu #1:#2\mskip1mu\bigr\}}
\newcommand{\sset}[1]{\bigl\{\mskip1mu#1\mskip1mu\bigr\}}
\def\deq{\ifmmode\equiv\else \macron\fi}
\newcommand{\defemph}[1]{{\sl #1}}
\let\temp=\colon\def\colon{\temp\mathopen{}}
\def\kG{{\mathcal G}}
\def\kQ{{\mathcal Q}}
\def\EE{{\mathbb E}}
\def\FF{{\mathbb F}}
\def\GG{{\mathbb G}}
\def\PP{{\mathbb P}}
\def\RR{{\mathbb R}}
\def\fg{{\mathfrak g}}
\newcommand{\bmat}[1]{\begin{bmatrix}#1\end{bmatrix}}
\renewcommand{\Im}{\operatorname{Im}}
\newcommand{\Shake}{{\sc\small SHAKE}}
\newcommand{\Rattle}{{\sc\small RATTLE}}
\renewcommand{\t}{{\mbox{\tiny $T$}}}

\ifx\pdftexversion\undefined
  \usepackage[dvips]{graphics}
\else
  \usepackage[pdftex]{graphics}
\fi
\numberwithin{equation}{section}

\newcommand{\ifpreprint}[2]{\ifthenelse{\equal\ispreprint{1}}{#1}{#2}}
\allowdisplaybreaks

\newtheoremstyle{slantbody}{\topsep}{\topsep}{\slshape\def\defemph{\defemphrm}}{}{\bfseries}{.}{6pt}{\thmname{#1}\thmnumber{ #2}\thmnote{ #3}}

\theoremstyle{slantbody}
\newtheorem{theorem}{Theorem}[section]

\newtheorem{lemma}[theorem]{\sf Lemma}

\begin{document}
\title{On converting any one-step method\\ to
  a variational integrator of the same order}

\author{
\llap{\parbox[t]{3.25in}{\it\small\begin{raggedleft}
George W.\ Patrick and Charles Cuell\\
Applied Mathematics and Mathematical Physics\\
Department of Mathematics and Statistics\\
\end{raggedleft}}}\qquad
\rlap{\parbox[t]{3.in}{\it\small\begin{raggedright}
Raymond J.\ Spiteri and William Zhang\\
Department of Computer Science\\\mbox{}\\
\end{raggedright}}}\\
\mbox{}\\[-10pt]
\it\small University of Saskatchewan, Saskatoon, Saskatchewan, S7N~5E6 Canada
}

\date{\it\small \today}

\maketitle

\vspace*{-.3in}\begin{abstract} 
\noindent In the formalism of constrained mechanics, such as that which
underlies the \Shake~and \Rattle~methods of molecular dynamics, we
present an algorithm to convert any one-step integration method to a
variational integrator of the same order. The one-step method is
arbitrary, and the conversion can be automated, resulting in a
powerful and flexible approach to the generation of novel variational
integrators with arbitrary order.

\end{abstract}

%
\section{\sf Introduction}
%
%
Consider a Lagrangian system defined by configuration space
$\hat\kQ\deq\sset{q}\deq\RR^N$, velocity phase space
$\!T\hat\kQ=\sset{(q,v)}=\RR^{2N}$, and Lagrangian
$L\colon\hat\kQ\to\RR$. Assume there is a holonomic constraint
$g(q)=0$, $g\colon\hat\kQ\to\RR^d$, suppose $g$ has full rank, and let
$\kQ=g^{-1}(0)$. The system evolves along curves $q(t)\in\kQ$ that
are critical points of the action
\[[eq:action]
  S\deq\int_a^b L\bigl(q(t),v(t)\bigr)\,dt
\]
subject to the fixed endpoint constraints $q(a),q(b)$ constant, and
the first order constraint $v(t)=q^\prime(t)$. This variational
principle is equivalent to the Euler--Lagrange equations
\[[eq:Euler-Lagrange-1]
  \frac{dq^i}{dt}=v^i,\quad\frac{dv^i}{dt}=A^i(q,v),
\]
where $A(q,v)$ is found by solving the linear (Lagrange multiplier)
problem (implicit sum on repeated indices)
\[[eq:Euler-Lagrange-2]
    &\frac{\partial^2L}{\partial v^i\partial v^j}A^j
    -\lambda_a\frac{\partial g^a}{\partial q^i}
    =-\frac{\partial^2L}{\partial v^i\partial q^j}v^j
    +\frac{\partial L}{\partial q^i},
  \\
    &-\frac{\partial g^a}{\partial q^i}A^i=
    \frac{\partial^2 g^a}{\partial q^i\partial q^j}v^iv^j.
\]
These are the general Lagrangian systems. For example, they
specialize to the Euler equations for the motion of a rigid
body~\cite{AbrahamR-MarsdenJE-1978-1, MarsdenJE-RatiuTS-1994-1} by
taking $\hat\kQ$ to be the $3\times3$ matrices $\sset{A}$ and $g^a$ to
the the upper triangular entries of $AA^\t-\!1$. And
they specialize to the Kirchhoff approximations for the motion of an
underwater 
vehicle~\cite{LeonardNE-MarsdenJE-1997-1,
              PatrickGW-RobertsRM-WulffC-2008-1}.
The same explicitly constrained formalism is exploited in the
molecular dynamics algorithms \Shake~and
\Rattle~\cite{AndersonHC-1983-1,
                    LeimkuhlerB-ReichS-2004-1,
                    LeimkuhlerBJ-SkeelRD-1994-1,
                    RyckaertJP-CiccottiG-BerendsenJC-1997-1}.

The objective here is symplectic and momentum-preserving simulation
of~\eqref{eq:Euler-Lagrange-1}. Such simulations may be systematically
generated by discretizing the defining variational principles, as
in~\cite{MarsdenJE-PatrickGW-ShkollerS-1998-1,
         MarsdenJE-WestM-2001-1,
         WendlandtJM-MarsdenJE-1997-1},
and particularly~\cite{LeyendeckerS-MarsdenJE-OrtizM-2008-1}, which
is specific to the constrained formalism.

The theory of variational integrators is elaborated 
in~\cite{CuellC-PatrickGW-2007-1, 
         CuellC-PatrickGW-2008-1},
based on geometric discretizations of the velocity phase space
$\!T\kQ$, i.e.,\ based on
\begin{enumerate}
  \item 
    certain assignments of \defemph{curve segments} to each
    $(q,v)\in\!T\hat\kQ$, where $(q,v)$ is tangent to $\kQ$; and
  \item
    a \defemph{discrete Lagrangian} $L_h$ that approximates $S_t$
    for $t=h$; $S_t$
    is defined to be the classical action~\eqref{eq:action} with $a=0$, $b=t$
    subject to~\eqref{eq:Euler-Lagrange-1}, and satisfies
  \[[eq:action-def]
    \frac{dS_t}{dt}=L(q^i,v^i).
  \]
\end{enumerate}
The discrete variational principle is a finite-dimensional constrained
optimization problem, in which the objective function is a sum of the
discrete Lagrangian on sequences in $\!T\kQ$. If the curve segments
associated with the elements of the sequence agree to order~$r$ with
the exact evolution of the Lagrangian system and the discrete
Lagrangian agrees to order~$r$ with the exact classical action, then
the variational integrator is order~$r$
accurate~\cite{PatrickGW-CuellC-2008-1}.

By definition, any one-step numerical integration method of order $r$
gives order~$r$ accurate solutions to
Equations~\eqref{eq:Euler-Lagrange-1} and~\eqref{eq:action-def}. Any
such method can be used to provide the curve segments and discrete
Lagrangian that are required to construct a variational integrator as
outlined above. In this article we derive the variational integrator
from the corresponding discrete Euler--Lagrange equations in terms of
a one-step integrator of ~\eqref{eq:Euler-Lagrange-1}
and~\eqref{eq:action-def}.

%
\section{Basic Algorithm}
%

Given a Lagrangian $L$, a constraint $g$, and a  one-step
numerical integrator of order $r$, which we call the \defemph{standard
  layer}, for the initial-value problem
\[
  &\frac{dq^i}{dt}=v^i,\quad\frac{dv^i}{dt}=A^i(q,v),\quad
    \frac{dS_t}{dt}=L(q,v),\\ 
  &q(0)=q,\quad v(0)=v,\quad S_t(0)=0,
\]
the aim is to generate a symplectic integrator of the same order. Let
the standard layer be represented by
\[
  (q,v)\mapsto R_t(q,v)=\bigl( R^q_t(q,v), R^v_t(q,v), 
  R_t^S(q,v)\bigr). \] Assume that the standard layer exactly
  preserves the constraint; this restriction will be lifted later. By
  differentiating $g(q)=0$, the space of vectors tangent to
  $\kQ=g^{-1}(0)$ is \[
  \!T\kQ\deq\set{w\deq(q,v)}{\!Dg(q)v=0}.
\]
Differentiating again,
\[[eq:tangent-constraint]
  \ifpreprint{}{&}
    \left.\frac d{d\epsilon}\right|_{\epsilon=0}\!Dg(q+\epsilon\,\delta q)
    (v+\epsilon\,\delta v) \ifpreprint{}{\\}
  \ifpreprint{}{&\quad\qquad\qquad}
    =v^\t\!D^2g(q)\,\delta q+\!Dg(q)\,\delta v,
\]
and so
\[
  \ifpreprint{}{&}
    \!T\!T\kQ\deq\set{\delta w\deq\bigl((q,v),(\delta q,\delta v)\bigr)}{
    \ifpreprint{}{\\&\quad\qquad}\!Dg(q)\,\delta q=0,\;  v^\t\!D^2g(q)\,\delta q+\!Dg(q)\,\delta v=0}.
\]
The notation $\!T_w\!T\kQ$ denotes the vector space $\set{\delta
w}{(w,\delta w)\in\!T\!T\kQ}$. $\!T\kQ$ and $\!T\!T\kQ$ are of
course the first and second tangent bundle of the constraint. The
notation for $\!D^2g(q)$ is problematic; it is a bilinear form with
values in $\RR^d$. The quantity $v^\t\!D^2g(q)$, where $v\in\RR^N$, denotes the
$d\times N$ matrix $v^i\partial^2g^a/\partial q^i\partial q^j(q)$.

We now construct the \defemph{symplectic layer} by defining the
following quantities.
\begin{enumerate}
\item 
  \emph{The \defemph{bias}:} a pair of numbers 
  \[
    -1\le\alpha^-\le0,\quad0\le\alpha^+\le1,
  \]
  such that $\alpha^+-\alpha^-=1$. 
\item
  \emph{The \defemph{time step}:} a number $h>0$.
\item 
  \emph{The maps $\partial^\pm_h$ and the discrete Lagrangian $L_h$:}
  the required curve segments are associated with each element $(q,v)$
  by $t\mapsto R^q_t(q,v)$ and the discrete Lagrangian is obtained
  from $(q,v)\mapsto R_h^S(q,v)$. $R^v_t$ is not used. The ends of the
  segments provide the maps
  \[
    &\partial^-_h(q,v)\equiv R^q_{h\alpha^-}(q,v),\quad
    \partial^+_h(q,v)\equiv R^q_{h\alpha^+}(q,v),
  \]
  and the discrete Lagrangian
  \[
    &L_h(q,v)\deq R^S_{h\alpha^+}(q,v)-R^S_{h\alpha^-}(q,v).
  \]
\item
  \emph{The time step of the symplectic layer:}
  Given $w_1=(q_1,v_1)\in\!T\kQ$,  solve the following \defemph{discrete
  Euler--Lagrange equations}~\cite{CuellC-PatrickGW-2008-1} for
  $w_2=(q_2,v_2)\in\!T\kQ$:
\[
  &\!DL_h(w_1)\,\delta w_1+\!DL_h(w_2)\,\delta w_2=0,\\
  &\partial^+_h(w_1)=\partial^-_h(w_2),
\]
for all $\delta w_1,\delta w_2$ satisfying
\[
  &\!D\partial_h^-(w_1)\,\delta w_1=0,\quad
    \!D\partial_h^+(w_2)\,\delta w_2=0,\\
  &\!D\partial_h^+(w_1)\,\delta w_1=\!D\partial_h^-(w_2)\,\delta w_2,\\
  &(w_1,\delta w_1)\in\!T\!T\kQ,\quad(w_2,\delta w_2)\in\!T\!T\kQ.
\]
\end{enumerate}
It is not necessary that the \emph{same} standard layer provides both
$\partial^-_h$ and $\partial^+_h$. For example, using any method to
construct $\partial^+_h$, the
\defemph{adjoint}~\cite{HairerE-LubichC-WannerG-2006-1} of the same
method to construct $\partial^-_h$, and the bias
$\alpha^-=\alpha^+=\frac12$, one evidently obtains a self-adjoint
symplectic layer. Self-adjoint methods respect time reversal in the
sense that a negative time step exactly reverses the discrete
evolution. Also, self-adjoint methods are necessarily of even order:
an odd-order self-adjoint method in fact has the next higher (even)
order of accuracy because its odd order truncation errors must equal
their negatives.

The symplectic layer corresponds to the finite-dimensional discrete
variational principle of finding the critical points of the
\defemph{discrete action}
\[
  S_h(w,\tilde w)\deq L_h(w)+L_h(\tilde w)
\]
subject to the constraints
\[
  &\partial_h^-(w)=\mbox{constant},\quad
    \partial_h^+(\tilde w)=\mbox{constant},\\
  &\partial_h^+(w)=\partial_h^-(\tilde w),\\
  &g(q,v)=0,\quad g(\tilde q,\tilde v)=0,
\]
as discussed in~\cite{CuellC-PatrickGW-2008-1}.

By the general theory, the maps $\partial_h^-,\partial_h^+$ split
$\!T\!T\kQ=\ker\!D\partial_h^-\oplus\ker\!D\partial_h^+$;
i.e.,\ $(w,\delta w)$ splits as
\[
  &\delta w=\delta w^++\delta w^-,\\
  &\delta w^+\in\ker\!D\partial_h^-(w),\quad
  \delta w^-\in\ker\!D\partial_h^+(w). \] (In the second line, the
  presence of plus with minus is intentional and conforms to the
  notation of~\cite{CuellC-PatrickGW-2008-1}.) The \defemph{discrete
    Lagrange one-form} is defined by \[
   \theta_{L_h}^-(w)\,\delta w\deq-\!DL_h(w)\,\delta w^-,
\]
and the general theory assures that the symplectic layer is a
symplectic integrator with respect to
$\omega_{L_h}\deq-\!d\theta_L^-$. 

In Lagrangian systems, the Noether theorem shows that the presence of
continuous symmetry is equivalent to the presence of conserved
momenta. For example, translational [rotational] symmetry implies
conservation of linear [angular] momentum;
see~\cite{AbrahamR-MarsdenJE-1978-1, MarsdenJE-RatiuTS-1994-1} for the
basic theory, some of which the discussion here must assume. The
discrete Noether theorem provides the same symmetry-momentum
equivalence for discrete Lagrangian systems: Suppose that a symmetry
group $\kG$ acts on $\hat\kQ$, such that
\begin{enumerate}
  \item $g$ is invariant;
  \item $\partial^+_h,\partial^-_h$ intertwine the lift of the action
  to $\!T\kQ$; and
  \item $L_h$ is invariant.
\end{enumerate}
Let $\fg$ be the Lie algebra of $\kG$, and let $\xi\in\fg$. Then the
symplectic layer preserves the \defemph{discrete momentum} defined by
\[
  J_\xi(w)\deq-\theta_{L_h}^-(w)\xi w,
\]
where $\xi\in\fg$ and $\xi w$ is the infinitesimal generator of $\xi$
at $w$. For example, $\kG$ could be a matrix Lie group that acts on
$\RR^N$ by matrix multiplication, the exact Lagrangian $L$ invariant,
the constraint $g$ invariant, and the standard layer an explicit
Runge--Kutta method. Then the standard layer intertwines the action on
$\!T\kQ$, and $L_h$ is invariant, and the symplectic layer will
preserve the corresponding discrete momenta.

Neither the discrete symplectic form nor the discrete momentum equals
in general the continuous counterpart. For example, if the system
is rotationally invariant then the symplectic layer need not preserve
some familiar mechanical angular momentum. The discrete Lagrangian
system has its own version of angular momentum, which is near to the
continuous one, but through its dependence on
$\partial_h^-,\partial_h^+$, is in general a complicated function of
$(q,v)$.

%
\section{Constrained algorithm}
%

We require the following standard lemma, which will justify the use of
a variety of (Lagrange) multipliers.

\begin{lemma}
Suppose that $A\colon\EE\to\FF$ and
$\mu\colon\EE\to\GG$ are linear. Let $\mu$ be onto. Then
$A(e)=0$ for all $e$ such that $\mu(e)=0$ if and only if there is a
$\lambda\colon\GG\to\FF$ such that $A=\lambda\mu$.
\end{lemma}

\begin{proof}
Suppose $A$ is zero on $\ker\mu$. Then $A$ drops to $\bar
A\colon\EE/\ker\mu\to\FF$. Also, $\mu$ drops to
$\bar\mu\colon\EE/\ker\mu\rightarrow\GG$, and this is a linear
isomorphism since $\mu$ is onto. Set $\lambda\equiv\bar
A\bar\mu^{-1}$. Then, letting $\pi\colon\EE\to\EE/\ker\mu$ be the
quotient map, $\lambda\mu(e)=\bar A\bar\mu^{-1}\mu(e)=\bar
A\pi(e)=A(e)$. Conversely, if $A=\lambda\mu$ and $\mu(e)=0$ then
$A(e)=\lambda\bigl(\mu(e)\bigr)=\lambda(0)=0$.
\end{proof}

We develop the equations used in the algorithm incrementally in
stages; see Figures~\ref{fig:levels-0-1-2} and~\ref{fig:levels-4}. In
the figures, the dimension counts for the equations and variables are
at right. The stages are equivalent representations of the same
algorithm, starting from the fundamental description of the algorithm
in Stage~0. Only the equations appearing in Stage~4 are implemented
and solved in practice. In going from Stage 0 to Stage 4, we increase
the number of equations to be solved from $5(N-d)$ to $12N+5d$.

\begin{figure*}[p]
  \footnotesize\begin{centering}
  \fbox{\parbox{.95\textwidth}{\vspace*{-8pt}
    \input{stage0123.tex}}}\end{centering}
  \caption{\label{fig:levels-0-1-2}\protect\ifpreprint{Stages 0--3.}{STAGES 0--3.}}
\end{figure*}

\begin{figure*}[p]
  \footnotesize\begin{centering}
  \fbox{\parbox{.95\textwidth}{\vspace*{-8pt}
  \input{stage4.tex}}}\end{centering}
  \caption{\label{fig:levels-4}\protect\ifpreprint{Stage 4: the derivatives in stage 3  are split into partial derivatives with respect to $q$ and $v$, and  the equations are reorganized. the approximates are correspondingly  the bottom primed equations.}{STAGE 4: THE DERIVATIVES IN STAGE 3  ARE SPLIT INTO PARTIAL DERIVATIVES WITH RESPECT TO $q$ AND $v$, AND  THE EQUATIONS ARE REORGANIZED. THE APPROXIMATES ARE CORRESPONDINGLY  THE BOTTOM PRIMED EQUATIONS.}}
\end{figure*}

\medskip\noindent\emph{Stage 0: The fundamental algorithm, given
  directly on the constrained phase space, as
  in~\cite{CuellC-PatrickGW-2008-1}.} The fundamental algorithm is defined
on $\kQ=g^{-1}(0)$, regarded as a submanifold of $\hat\kQ$. The
constraints \[
  \!D\partial_h^-(w_1)\,\delta w_1=0,\quad
  \!D\partial_h^+(w_2)\,\delta w_2=0,
\]
are enforced with  multipliers (in $\!T\kQ$) $\lambda^-$ and
$\lambda^+$, both having dimension $\dim\kQ=N-d$. The constraint
\[
  \!D\partial_h^+(w_1)\,\delta w_1=\!D\partial_h^-(w_2)\,\delta w_2,\
\]
is enforced with a multiplier $\mu$. Setting $\delta w_1$ and $\delta
w_2$ alternately to zero gives the two equations~(S0.2) and~(S0.1),
respectively. The \defemph{connecting constraint}
$\partial^+_h(w_1)=\partial^-_h(w_2)$ translates unchanged to (S0.3).

\medskip\noindent\emph{Stage 1: Enforce the restriction $\delta
  w_i\in\!T_{w_i}\!T\kQ$ by introducing multipliers.} The restriction
$(\delta q,\delta v)=\delta w\in\!T_w\!T\kQ$ is \[
  \!Dg(q)\,\delta q=0,\quad 
  v^\t\!D^2g(q)\,\delta q+\!Dg(q)\,\delta v=0.
\]
For each equation in~(S0.1) and (S0.2), there are two corresponding
multipliers $\nu_1,\nu_2$, which are row vectors of length $d$.

\medskip\noindent\emph{Stage 2: Disambiguate
  $\lambda^-,\lambda^+,\mu$.} The multipliers
$\lambda^+,\lambda^-,\mu$ are ambiguous as row vectors in $\RR^N$ up
to any vector orthogonal to the constraint. Specifying  them 
into $\!T\kQ$ disambiguates them. Because $\kQ=g^{-1}(0)$ and the rows
of $\!Dg$ span the orthogonal complement to the tangent space of
$\kQ$, the multipliers should have zero dot product with the rows of
$\!Dg$, e.g.,\ $\lambda\!Dg(q)^\t=0$ for a multiplier $\lambda$ at
$q\in\kQ$.

\medskip\noindent\emph{Stage 3: Lift the restriction that the standard
  layer preserves the constraint.} A numerical integrator of order $r$
used in the standard layer will in general preserve the constraint $g$
only to accuracy order~$r$. We posit a map
$\iota\colon\RR^N\times\RR^d\to\RR^N$ such that \begin{enumerate}
\item $\iota(q,0)=q$; and \item if $g(q)=0$, then
  $\Im\!D_\theta\iota(q,0)$ is a complement of $\ker\!Dg$; i.e.,\
  $\Im\!D_\theta\iota(q,0)\oplus\ker\!Dg(q)=\RR^N$ for all $q\in\kQ$.
\end{enumerate} 
Define a map $\PP$ by
\[[eq:def-general-project]
  \PP(\hat q)\deq q,\quad
  \hat q=\iota(q,\theta),\quad g(q)=0;
\]
see \ifpreprint{the below figure.}{Figure~\ref{fig:fig1}.}  The map $\PP$ follows the
constant~$\theta$-fibers of $\iota$ to where they intersect with
$\kQ$.  The maps $\partial_h^+,\partial_h^-$ are defined by the
(constraint preserving) standard layer above. Letting the
unconstrained standard layer define
$\hat\partial_h^+,\hat\partial_h^-$, we redefine
$\partial_h^+,\partial_h^-$ by
\[
  \partial_h^+\deq\PP\hat\partial_h^+,\quad
  \partial_h^-\deq\PP\hat\partial_h^-.
\]
Also, we introduce the new variables
\[
  &q_1^-\deq\partial_h^-(w_1),\quad
  q_2^+\deq\partial_h^+(w_2),\\
  &\bar q\deq\partial_h^+(w_1)=\partial_h^-(w_2),
\]
with variable and equation count each of $N$.

\begin{figure}[hb]
  \begin{center}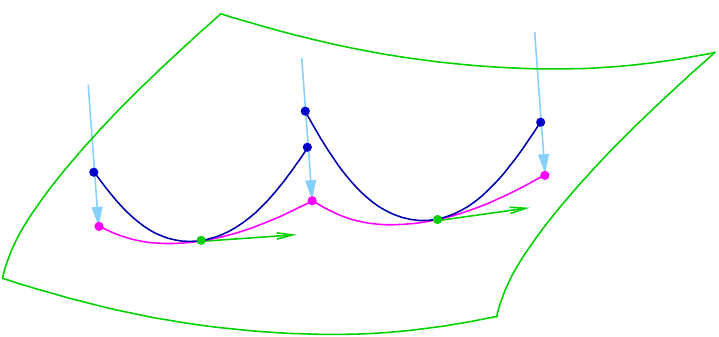\end{center}
  \ifpreprint{}{\caption{\label{fig:fig1}}}
\end{figure}

Usually an explicit projection $\PP$ is not available and has to be
computed iteratively. For example, $\iota$ and $\PP$ may be naturally
defined by
\[[eq:def-project]
  &\PP(\hat q)=q,\\
  &\hat q=q+\!Dg\bigl(q\bigr)^\t\theta,\quad\mbox{$\theta\in\RR^d$},\\
  &g(q)=0.
\]
As $\theta$ varies, this particular $\iota(q,\theta)$ moves $q\in\kQ$
away from $\kQ$ orthogonally; the reverse, obtained from $\PP$,
projects $\hat q$ orthogonally to $\kQ$.

The connecting equation~(S2.3) must be assumed to be full rank into
$\kQ$; hence its (linearly independent) equation count is $N-d$. There
are actually $N$~equations when the image is considered into $\RR^N$,
as it must be for computations, but $d$ of those are redundant because
$w_2\in\kQ$ by assumption and $\partial^+_h,\partial^-_h$ preserve the
constraints. In the second of (S3.5), explicitly writing the
projection using new variable $\theta^+$ resolves this problem because
$\theta^+$ robustly moves~$q$ away from the constraint~$\kQ$, i.e.,\
locally linearly in a nondegenerate way. So Equation~(S3.5) replaces
(half of) Equation~(S2.3), with equation count $N$ rather than $N-d$,
while the number of variables increases by $d$ because that is the
count for $\theta^+$.

\medskip\noindent\emph{Stage 4: Split (S3.1) and (S3.2) into partial
  derivatives with
  respect to $q$ and $v$; rearrange terms and group equations.}
The $q$ and $v$ partial derivatives of Equation~(S3.1) give
Equations~(S4.3a) and~(S4.4a), respectively. Similarly, partial
derivatives of~(S3.2) give~(S4.5a) and~(S4.6a). We group the equations
so they can (eventually) approximated by linear equations with the
same coefficient matrices, as will be seen.

%
\section{Implementation}
%

We present here a strategy for the Stage~4 computation.

\subsection{The vector field and its derivatives.} We now specialize
to Lagrangians of the form \[[eq:lagform]
L\deq\frac12m_{ij}(q)v^iv^j+a_j(q)v^j-V(q). \] This is the most
general quadratic Lagrangian with configuration-dependent
coefficients. As is easily verified,
Equation~\eqref{eq:Euler-Lagrange-2} becomes \[[eq:impl]
  &m_{ij}A^j-\lambda_a\frac{\partial g^a}{\partial q^i}
    =-\Gamma_{ikl}v^kv^l-b_{ij}v^j-\frac{\partial V}{\partial q^i},\\
  &-\frac{\partial g^a}{\partial q^i}A^i
    =\frac{\partial^2 g^a}{\partial q^i\partial q^j}v^iv^j,\\
  &\Gamma_{ikl}\deq\frac12\left(\frac{\partial m_{il}}{\partial q^k}
    +\frac{\partial m_{ik}}{\partial q^l}
    -\frac{\partial m_{kl}}{\partial q^i}\right),\\
  &b_{ij}\deq\frac{\partial a_i}{\partial q^j}
    -\frac{\partial a_j}{\partial q^i}.\\
\]
These are all linear equations for $A$ and $\lambda$ with 
coefficient matrix of the  form
\[[eq:KKT1]
  \bmat{M(q)&-\!Dg(q)^\t\\-\!Dg(q)&0},\quad M(q)\deq\bigl[m_{ij}(q)\bigr].
\]

The algorithm requires the derivatives of the maps $\partial^-_h$ and
$\partial^+_h$ as well as the derivative of $L_h$. Automatic
differentiation~\cite{GriewankA-2000-1} can be used to compute these
by computing the derivative of the one-step method in the standard
layer that defines them. Alternatively, by Lemma~4.1
of~\cite{HairerE-LubichC-WannerG-2006-1}, the derivative of the
standard layer $R_t$ that uses a Runge--Kutta method may be computed
as the same Runge--Kutta method applied to the equations of first
variation. In this case, it is only required to determine the
derivative of the vector field.

To compute the derivative of $A^j$ with respect to $q^m$ and $v^m$, we
differentiate the first equation in~\eqref{eq:impl}:
\[
  &
  m_{ij}\frac{\partial A^j}{\partial q^m}
    -\frac{\partial\lambda_a}{\partial q^m}\frac{\partial g^a}{\partial q^i}
  \ifpreprint{}{\\}
  \ifpreprint{}{&\qquad}=-\frac{\partial \Gamma_{ikl}}{\partial q^m}v^kv^l
    -\frac{\partial b_{ij}}{\partial q^m}v^j
    -\frac{\partial^2V}{\partial q^i\partial q^m}
  \ifpreprint{}{\\&\qquad\qquad}
    -\frac{\partial m_{ij}}{\partial q^m}A^j
    -\frac{\partial^2g^a}{\partial q^i\partial q^m}\lambda_a,\\
  &m_{ij}\frac{\partial A^j}{\partial v^m}
    -\frac{\partial\lambda_a}{\partial v^m}\frac{\partial g^a}{\partial q^i}
    =-2\Gamma_{ikm}v^k
\]
and for the constraints
\[
  &-\frac{\partial g^a}{\partial q^i}\frac{\partial A^i}{\partial q^m}=
    \frac{\partial^3g^a}{\partial q^i\partial q^j\partial q^m}v^iv^j
    \frac{\partial^2g^a}{\partial q^i\partial q^m}A^i,\\
  &
    -\frac{\partial g^a}{\partial q^i}\frac{\partial A^i}{\partial v^m}=
    2\frac{\partial^2g^a}{\partial q^i\partial q^m}v^i.
\]
These resulting equations are linear in the required derivatives with
coefficient matrix~\eqref{eq:KKT1}.

\subsection{Fixed-point iteration} A possible approach for solving the
Stage~4 (implicit) equations is by a fixed-point iteration. To find
solutions to a generic equation $f(x)=0$, split $f=f_0-\Delta f$
such that, for all $b$, an explicit solution to the equation
$f_0(x)=b$ is available. We call $f_0$ an \defemph{approximate}. If
$\Delta f$ is sufficiently small, then in a
suitable neighbourhood of the solution, the iteration
$x_{i+1}=f_0^{-1}\bigl(\Delta f(x_i)\bigr)$ converges to a solution of
$f(x)=0$: \[
  f_0(x)=\Delta f(x)=f_0(x)-f(x)\quad\Leftrightarrow\quad f(x)=0.
\]
The iteration is, given an initial iterate $x_0$,
\[
  &r_{i+1}=f_0(x_i)-f(x_i),\quad \mbox{solve $f_0(x_{i+1})=r_{i+1}$},
\]
or equivalently, after substituting $f_0(x_i)=r_i$,
\[[eq:generic-iteration]
  &r_0=f_0(x_0),\\
  &r_{i+1}=r_i-f(x_i),\\
  &\mbox{solve $f_0(x_{i+1})=r_{i+1}$}.
\]
This approach is useful for the Stage~4 equations because they are
nonlinear and a good choice for $f_0$ is generally available. In this
way, the Stage~4 computation may be organized into its equations and
corresponding approximates. The required solution is obtained by
iteratively evaluating the equations themselves and then solving for
the approximates using~\eqref{eq:generic-iteration}.

\subsection{Stage 4 approximates}

Because the time step is small, the various configurations $q_1,\bar
q,q_2$, etc., are all close. Let $G_0$, $M_0$, and $a_0$, be
approximations to $\!Dg(q)$, $M(q)$, and $a(q)$ respectively, obtained
by evaluation at some such configuration; e.g., the configuration
$q_1$ is a likely candidate.

Equations~(S4.1), (S4.2), and~(S4.8) are of the form of
Equations~\eqref{eq:def-project}, which, for sufficiently small $h$,
can be effectively approximated by Taylor expansion of $g(q)=0$ at
$q=\hat q$:
\[
  &(\hat q-q)-G_0^\t\theta=0,\\
  &-G_0(\hat q- q)-g(\hat q)=0.
\]
Equations (S4.3c), (S4.4c), (S4.5c), and (S4.6c) involve the derivative of
$\PP$ in expressions such as
\[[eq:tilde-lambda]
  \hat\lambda^\t=\!D\PP(q)\lambda^\t.
\]
Differentiating equations~\eqref{eq:def-project} gives
\[[eq:diff-project]
  &\delta\hat q=\delta q+\!Dg(q)^\t\,\delta\theta
  +\!D^2(\theta^\t g)(q)\,\delta q,\\
  &\delta g=\!Dg(q)\,\delta q.
\]
The matrix $\!D\PP(q)$ is obtained by discarding $\delta\theta$ after
the inverse of \eqref{eq:diff-project}, with $\delta g=0$, i.e.,\
\[
  \!D\PP(q)=
  \bmat{\!1&0}\bmat{\!1+\!D^2(\theta^\t g)(q)&\!Dg(q)^\t\\\!Dg(q)&0}^{-1}
  \bmat{\!1&0\\0&0},
\]
and Equation~\eqref{eq:tilde-lambda} becomes
\[
  \bmat{\hat\lambda^\t\\0}&=
  \bmat{\!1&0\\0&0}
  \bmat{\!1+\!D^2(\theta^\t g)(q)&\!Dg(q)^\t\\\!Dg(q)&0}^{-1}
\bmat{\lambda^\t\\0}.
\]
If we define a variable $z$ and put it in place of the zero in the
matrix at left, then we can put a unit matrix in the $(2,2)$ slot of
the first matrix on the right side, and invert. The result is the
linear equation
\[[eq:lambda-pm-mu-project]
  \bmat{\!1+\!D^2(\theta^\t g)(q)&\!Dg(q)^\t\\\!Dg(q)&0}\bmat{\hat\lambda^\t\\z}=
  \bmat{\lambda^\t\\0}.
\]
This replaces (S4.3c), (S4.4c), (S4.5c), and (S4.6c), and one can use the
approximate
\[[eq:KKT-1]
  \bmat{\!1&G_0^\t\\G_0&0}
  \bmat{\hat\lambda^\t\\z}=\bmat{\lambda^\t\\0},
\]
which is computationally equivalent to
\[
  \bmat{\!1&-G_0^\t\\-G_0&0}
  \bmat{-\hat\lambda^\t\\z}=\bmat{-\lambda^\t\\0}.
\]

The remaining approximates are driven by the basic data of the
variational principle:
\[[eq:approx-system]
  &\!D_qL_h\approx h\!D_qL,\quad
  \!D_vL_h\approx hv^\t M_0+ha_0,\\
  &\hat\partial_h^+(q,v)\approx q+h\alpha^+v,\quad
  \hat\partial_h^-(q,v)\approx q+h\alpha^-v.
\]
In order, the approximates for Equations~(4.3a,b) obtained from
\eqref{eq:approx-system}, $\hat\lambda^-\approx\lambda^-$, and
$\hat\mu\approx\mu$ are as follows:
\[
  &\!D_q\partial^-_h(w_1)\approx\!1,\quad
     \hat\lambda^-\!D_q\partial^-_h(w_1)\approx\lambda^-,\\
  &\!D_q\partial^+_h(w_1)\approx\!1,\quad
     \hat\mu_1\!D_q\partial^+_h(w_1)\approx\mu,\\
  &\nu^-_1\!Dg(q_1)\approx\nu^-_1 G_0=0,\\
\]
resulting in
\[
  &\lambda^-+\mu+\nu_1^-G_0+\nu_2^-v_1^\t\!D^2g(q_1)
  -h\!D_qL(q_1,v_1)=0,\\
  &\lambda^-G_0^\t=0. \] These are equations (S4.3a$^\prime$) and
  (S4.3b$^\prime$). Similarly one obtains the approximates
  (S4.4a$^\prime$b$^\prime$)--(S4.6a$^\prime$b$^\prime$), noting
  however that
  in~(S4.5a$^\prime$) there are the further approximations 
  \[
    \nu_2^+v_1^\t\!D^2g(q_1)\approx\nu_2^-v_1^\t\!D^2g(q_1),
    \ifpreprint{\qquad}{\quad}
    \!D_qL(q_2,v_2)\approx\!D_qL(q_1,v_1).
\]

Equations (S4.7a) and (S4.7b) are
\[
\hat\partial^-_h(w_2)-\!Dg(\bar q)^\t\theta^+=\bar q,\quad
g(q_2)=0,
\]
which have to be solved for $q_2$ and $\theta^+$. Taylor expanding the
second equation at $\bar q$, and using $g(\bar q)=0$, gives the approximate
\[
&q_2-\bar q-\!Dg(q_1)^\t\theta^++h\alpha^+v_2=0,\\
&G_0(q_2-\bar q)=0,
\]
which are Equations~(S4.7a$^\prime$b$^\prime$).

\subsection{Stage 4 solution}

The approximates~(S4.3a$^\prime$b$^\prime$)
through~(S4.5a$^\prime$b$^\prime$) are solvable for the multipliers
\[
  \lambda^-,\;\mu,\;\lambda^+,\;\nu_1^-,\;\nu_2^-,\;\nu_1^+,
\]
and the linear equations are all of the form~\eqref{eq:KKT-1}. Indeed,
one adds $\alpha^-$ times~(S4.3b$^\prime$) and $\alpha^+$
times (S4.4b$^\prime$), and then that together with~(S4.4a$^\prime$)
can be solved for $\alpha^-\lambda^-+\alpha^+\mu$ and $\nu_2^-$.
Then, and similarly, (S4.3a$^\prime$b$^\prime$) can be solved for
$\lambda^-+\mu$ and $\nu_1^-$. Together these give
$\lambda^-,\mu,\nu_1^-,\nu_2^-$. Because $\mu$ is then known,
(S4.5a$^\prime$) with (S4.5b$^\prime$) minus
(S4.4b$^\prime$) provide $\lambda^+$ and $\nu_1^+$. 
In the same way~(S4.6a$^\prime$b$^\prime$) may be used to update
$\nu_2^+$ and $v_2$. Finally, (S4.7a$^\prime$b$^\prime$) are solved
for $q_2$ and
\[
  \hat\lambda^-,\;\hat\mu_1,\;\hat\mu_2,\;\hat\lambda^+
\]
occurring in~(S4.3c), (S4.4c), (S4.5c), and (S4.6c$^\prime$), may all
be updated using appropriate versions
of \eqref{eq:lambda-pm-mu-project} and its approximates. 
The entire procedure can then be iterated until the variables $q_2$
and $v_2$ are at a predetermined accuracy.

If should be noted that the multipliers that impose the constraints
$g$, i.e.,\ $\nu_1^-,\nu_2^-,\nu_1^+,\nu_2^+$ are not unique. For
example, doubling $g$ results in halving these multipliers. Such
multipliers are nonphysical, and convergence of the iteration of
Stage~4 should not be bound to the convergence of the multipliers
themselves. Rather, the degree of convergence can be determined from
products such as $\nu_2^-\!Dg(q_1)$, which generally have the physical
meaning
of \defemph{force of constraint}; i.e.,\ they are added directly in the
equations to quantities with a physical interpretation.

\bibliographystyle{plain}\footnotesize
\bibliography{0references/abbrev,0references/ref}
\end{document}